\documentclass[11pt,a4paper]{amsart}

\usepackage{amsfonts}
\usepackage{amscd}
\usepackage{amssymb}
\usepackage{amsmath}
\usepackage{latexsym}
\usepackage[dvips]{graphicx}
\usepackage{rotating}
\usepackage{tabularx}

\def\Conv {{ \mbox{\rm Conv}}}

\def\Vol {{\mbox{\rm Vol}}}

\def\init {{\rm \mbox{init}}}
\def\Res {\mbox{\rm Res}}

\def\MV {\mbox{\rm MV}}

\def\C{\mathbb{C}}

\def \N{\mathbb{N}}
\def \P{\mathbb{P}}
\def \Q{\mathbb{Q}}
\def \R{\mathbb{R}}
\def \Z{\mathbb{Z}}

\def\cA {{\mathcal A}}
\def\cB {{\mathcal B}}

\def\cE {{\mathcal E}}

\def\cL {{\mathcal L}}
\def\cM {{\mathcal M}}

\newtheorem{lem}{Lemma}[section]

\newtheorem{thm}[lem]{Theorem}

\newtheorem{defn}{Definition}[section]

\newtheorem{rem}[defn]{Remark}
\newtheorem{exmpl}[defn]{Example}

           {\vspace{3.3mm}
           \noindent{\bf #1}\it}%
           {\vspace{3.3mm}}

\newenvironment{Proof}[1]{{\it #1}}{\hfill\mbox{$\Box$} }

\begin{document}

{\title[The height of the mixed  sparse resultant]{The 
height of the mixed sparse resultant}
\author{Mart{\'\i}n Sombra}
  
\address{Universit{\'e} de Paris 7, UFR de Math{\'e}matiques,
{\'E}quipe de G{\'e}om{\'e}trie et Dynamique, 
2 place Jussieu, 75251 Paris Cedex 05,
France; 
and
Departamento de Matem{\'a}tica,
Universidad Nacional de La Plata,
Calle 50 y 115,
1900 La Plata, Argentina.
\medskip  } 

\email{sombra@math.jussieu.fr}

\date{November 28, 2002} 



\thanks{Supported by a Marie Curie Post-doctoral fellowship 
of the European Community Program 
{\em Improving  Human
Research Potential and the Socio-economic Knowledge Base}, 
contract n\textordmasculine \ HPMFCT-2000-00709.} 

\abstract{We present an upper bound  for the height 
of the 
mixed sparse resultant, 
defined as the logarithm of the maximum
modulus of its coefficients. 
We obtain a similar estimate for its Mahler measure.}}

\maketitle

\vspace{-8mm}


\setcounter{section}{1}

\section*{}

Let $\cA_0, \dots, \cA_n \subset \Z^n $ be finite sets 
of integer vectors and let
\,$ 
\Res_{\cA_0, \dots, \cA_n} \in \Z[U_0, \dots, U_n ] 
$\, 
be the associated mixed sparse resultant 
--- or $ (\cA_0, \dots, \cA_n)$-resultant --- which is a 
polynomial in $n+1$ groups
\,$U_i := \{ U_{i \, a} \, ; \, a \in \cA_i \} $\, of 
\,$m_i := \# \cA_i$\, variables each. 
We refer to \cite{Sturmfels94} 
and \cite[Chapter 7]{CoLiOS98}
for the definitions and basic facts. 

This resultant is widely used as a  tool 
for 
polynomial equation
solving, 
a fact that sparked a lot interest in its  computation,
see e.g. \cite[Sec. 7.6]{CoLiOS98}, \cite{EmMo99}, \cite{DAndrea01}, 
\cite{JeKrSaSo02},  
while  
it is also studied from a more  theoretical point of view 
because of its connections with toric
varieties
and hypergeometric functions, 
see e.g. \cite{GeKaZe94}, \cite{CaDiSt98}.

\bigskip 

We assume for the sequel that the family of supports 
 $\cA_0, \dots, \cA_n$ is {essential} 
(see \cite[Sec. 1]{Sturmfels94})
which  
does not represent any loss
of generality, by \cite[Cor. 1.1]{Sturmfels94}.  

\smallskip 

Set $\cA:= (\cA_0, \dots, \cA_n)$, and let $L_\cA \subset \Z^n$ 
denote the $\Z$-module affinely spanned by the pointwise  
sum $\sum_{i=0}^n \cA_i$. 
This is  a 
subgroup of $\Z^n$ of finite 
index  
$$ 
[\Z^n : L_\cA] := \# (\Z^n/L_\cA) 
$$
because we assumed that the family $\cA$ is essential. 
Also set 
$
Q_i := \Conv(\cA_i) \subset \R^n$ 
for the convex hull of $\cA_i$ 
for $i=0, \dots n$. 

We note by $\MV$ the {mixed volume} function 
as defined in e.g. 
\cite[Sec. 7.4]{CoLiOS98}: 
this 
is normalized so that for a 
polytope $P \subset \R^n$, 
the mixed volume $\MV(P, \dots, P)$ equals $n!$ times the Euclidean  
volume  $ \Vol_{\R^n} (P)$. 
We also  set $\Vol(P):= \MV(P, \dots, P)=  n!\, 
\Vol_{\R^n} (P)$.

\smallskip 

Under this 
  notation and assumption, the resultant 
is a multihomogeneous polynomial of degree
$$
\deg_{U_i}\hspace{-1mm}\Big(\Res_{\cA_0, \dots, \cA_n} \Big) \, = \,  
\frac{1}{[\Z^n : L_\cA] } \, \MV (Q_0, \dots, 
Q_{i-1} , Q_{i+1} , \dots, Q_n) \ >0 
$$
with respect to each group of variables $U_i$, see 
\cite[Cor. 2.4]{PeSt93}. 

\smallskip

The {\em absolute height} of a polynomial $g =\sum_{a} c_a \, x^a
\in \C[x_1, \dots, x_n]$ 
is defined as 
$H(g):= \max \{ |c_a| \, ; \, a \in \N^n \}$. 
Hereby we will be mainly concerned with its 
{\em (logarithmic)  height}: 
$$
h(g) := \log H(g)  =
\log \max \{ |c_a| \, ; \, a \in \N^n \}. 
$$

\smallskip 

\smallskip 

The main result of this paper is the following upper bound for the 
height of the resultant:

\begin{thm} \label{h}
$$ 
h \Big(\Res_{\cA_0, \dots, \cA_n} \Big) \, \le \, 
\frac{1}{[\Z^n:L_\cA]} 
\, \sum_{i=0}^n \MV(Q_0, \dots, Q_{i-1}, Q_{i+1} , \dots, Q_n) \, 
\log (\# \cA_i). 
$$ 
\end{thm} 

We write for short  \,$\Res_\cA:= \Res_{\cA_0, \dots, \cA_n}$\,
and 
\,$
MV_i(\cA):= 
\frac{1}{[\Z^n:L_\cA]} \, \MV(Q_0, \dots, Q_{i-1}, Q_{i+1} , $ 
$\dots, Q_n)  
$\,
for $i=0, \dots, n$. 
Thus the previous result can be rephrased as 
$$
H\Big(\Res_{\cA} \Big) \, \le \, 
\prod_{i=0}^n 
(\# \cA_i)^{MV_i(\cA)}.  
$$

\smallskip 

This improves our previous bound for the
unmixed case \cite[Cor. 2.5]{Sombra02}
and  extends it to the general case. 
We remark that the obtained  
upper bound is {\em polynomial} in the size of the input family
of supports
\,$\cA$\, and  in the  mixed volumes \,$MV_i(\cA) $\,, 
 and hence it represents a 
truly substantial improvement over all previous
general estimates. 
These are 
the ones which follow either  from the Canny-Emiris type
formulas 
(see 
 Inequality (\ref{c-e}) in the 
appendix, see also  
\cite[Prop. 1.7]{KrPaSo01}, \cite[Thm. 23]{Rojas00})
or from direct application of the unmixed case 
(see the inequality (\ref{unmixed}) below for $k=1$).

\medskip 

We also consider the 
Mahler measure, which is another  usual notion for the size of a
$n$-variate polynomial. 
The {\em Mahler measure} of  
$g\in \C[x_1, \ldots, x_{n}]$ 
is  defined as 
$$
m(g) 
:= 
\int_{S_1^n}  \log |g| \ d\mu^n , 
$$ 
where  $S_1 \subset \C$ 
is the unit circle, and $d\mu$ is the Haar measure over $S_1$ 
of total mass $1$. 
This can be compared with the height: 
in our case
\begin{equation} \label{mh} 
-\sum_{i=0}^n MV_i(\cA) \,  \log (m_i)  
\, \le \, 
m(\Res_\cA) - h(\Res_\cA) 
\, \le \,  
\sum_{i=0}^n MV_i(\cA) \,  \log (m_i) 
\end{equation} 
by \cite[Lem. 1.1]{KrPaSo01}. 
We refer to \cite[Sec. 1.1.1]{KrPaSo01} for an account on some of the 
notions of height of complex  polynomials: just note that the height  
$h(g)$ here coincides with  $\log|g|_\infty$ in that  
reference.

\smallskip

We obtain the same estimate as
before for the Mahler measure of the resultant.

\begin{thm} \label{m} 
$$
m \Big( \Res_{\cA_0, \dots, \cA_n} \Big)  
\, \le \, 
\frac{1}{[\Z^n:L_\cA]} 
\, \sum_{i=0}^n \MV(Q_0, \dots, Q_{i-1}, Q_{i+1} , \dots, Q_n) \, 
\log (\# \cA_i).
$$ 
\end{thm} 

Note that this improves by a factor of 2 the estimate which would derive 
from direct application of Theorem~\ref{h} and the inequalities 
(\ref{mh}) above. 

\medskip 

Both estimates are a consequence of the following:

\begin{lem} \label{1}

Let $f_0 \in \C^{\cA_0}, \dots, f_n\in \C^{\cA_n} $. 
Then 
$$
\log \Big| \Res_{\cA_0, \dots, \cA_n} (f_0, \dots, f_n)  \Big| 
\, \le \, \frac{1}{[\Z^n:L_\cA]} 
\,  \sum_{i=0}^n \MV(Q_0, \dots, Q_{i-1}, Q_{i+1} , \dots, Q_n) 
\, \log ||f_i||_1,
$$ 
where $||f_i||_1 := \sum_{a \in \cA_i} |f_{i \, a} | $ denotes the
$\ell_1$-norm of the vector $f_i = (f_{i\, a} \, ; \, a \in \cA_i)$.

\end{lem} 

Let
$g =\sum_{a} c_a \, x^a \in \C[x_1, \dots,
x_n]$.  
Then for $a\in \N^n$ we have that  
$$
c_a = \int_{S_1^n} 
\frac{g(z_1, \dots, z_n)}{z_1^{a_1+1} \cdots z_n^{a_n+1} } 
\ d\mu^n 
$$
by  Cauchy's formula and so 
\,$ h(g) \le \sup \Big\{ \log |g(\xi)| \, ; \, \xi \in S_1^n
\Big\}$. 
Thus Theorem~\ref{h}
is a consequence of this inequality applied to $g:= \Res_\cA$, 
together with Lemma~\ref{1}.

On the other hand, Theorem~\ref{m} follows from
Lemma~\ref{1} 
by a straightforward estimation of the  integral in the definition of the 
Mahler measure.

\medskip

\begin{Proof}{Proof of Lemma~\ref{1}.--} 
Let $ k \in \N$. 
Then let $k \, \cA_i \subset \Z^n$ 
denote the  pointwise sum of $k$ copies of $\cA_i$, and  
set $k\, \cA:= (k \, \cA_0, \dots, k \, \cA_n) $. 
It is easy to verify that $k\, \cA$ is also essential,  $L_{k \,\cA} =
L_\cA$ 
and 
$\Conv(k \, \cA_i) = k \, Q_i$.

We identify each  $f_i \in \C^{\cA_i} $ with the corresponding  Laurent  
polynomial 
$
f_i= \sum_{a \in \cA_i} f_{i\,a} \, x^a$, and we 
set 
$f_i^k \in \C^{k  \cA_i } $ for the vector which corresponds to the
$k$-th power of $f_i$.  
By the factorization formula for resultants  
\cite[Prop. 7.1]{PeSt93} we get that   
$$
\Res_{k\cA} (f_0^k, \dots, f_n^k) 
=
\Res_\cA(f_0, \dots, f_n)^{k^{n+1}}
$$
and so 
\begin{eqnarray} \label{prod}  
k^{n+1} \, \log\Big|\Res_\cA(f_0, \dots, f_n)\Big|  
& \le &  
h(\Res_{k  \cA}) +   
\sum_{i=0}^n
MV_i(k\, \cA) \, \log ||f_i^k||_1  \nonumber \\[-2mm]  
& \le &   
h(\Res_{k \cA}) + 
k^{n+1} \, \sum_{i=0}^n MV_i(\cA)  \, 
 \, \log ||f_i||_1. 
\end{eqnarray} 
The first inequality follows from the straightforward estimate
$|G(u_0, \dots, u_n) | \le H(G) \, 
$ \newline $ \prod_{i=0}^n ||u_i||_1^{d_i} $
for a multihomogeneous polynomial $G$ of degree $d_i$ in each group of 
variables, 
applied to $G:= \Res_{k \, \cA} $ and $u_i := f_i^k$. 
The second one follows  from the linearity of the mixed volume, 
and the sub-additivity of the $\ell_1$-norm with 
respect to polynomial multiplication 
(which implies that $\log||f_i^k||_1 \le k\, \log
||f_i||_1$). 

\medskip 

Now let \,$\cB \subset \Z^n$\, be any finite set such
that $L_\cB= \Z^n$ and 
such that $\cA_0, \dots,
\cA_n \subset \cB $.
Set 
\,$n(k):= \# k\, \cB$\, and  $ P:= \Conv(\cB)  \subset \R^n$. 
Then the (unmixed) resultant \,$\Res_{k\cB} $\, is a polynomial in 
 \,$(n+1) \, n(k) $\, variables and total degree 
\,$(n+1) \,\Vol(k\, P )= 
(n+1) \, k^n\,\Vol (P)$.
We have also that  $L_{k \cB} = \Z^n$ 
and so we are in the hypothesis of 
\cite[Cor. 2.5]{Sombra02}, which gives the height estimate 
$$ 
h(\Res_{k\cB})  \le   
2 \, (n+1)  \,  \log\hspace{-1mm}\Big(n(k) \Big) \, 
\Vol (k\, P)
=
2 \, (n+1)  \,  \log\hspace{-1mm}\Big(n(k) \Big) \, 
 k^n\,\Vol (P). 
$$
We have that 
$k \, \cA_i \subset k\, \cB $ for $i=0, \dots, n$ and so 
by \cite[Cor. 4.2]{Sturmfels94}
 there exists a monomial order $\prec$ 
such that 
$ \Res_{k\cA}  $ 
divides the initial form \,$ \init_\prec (\Res_{k \, \cB} ) $. 
This is a polynomial in \,$(n+1) \, n(k)$\, variables of degree and height
bounded by those of $\Res_{k\cB}$, and so 
\begin{eqnarray} \label{unmixed} 
h(\Res_{k \cA}) & \le &  
h(\Res_{k\cB} ) + 2\, 
\log\Big((n+1) \, n(k)  +1\Big) 
\, (n+1)  \,  k^n \, 
\Vol (P) \nonumber \\[2mm] 
& \le &     4\, (n+1)
\, \log \Big( (n+1) \, n(k)  +1 \Big) \, k^n\,\Vol (P)
\end{eqnarray} 
by the inequality  $h(f) \le h(g) + 2 \, \deg(g) \, \log(n+1)$, which 
holds for 
$f,g \in \Z[x_1, \dots, x_n]$ such that $f|g$ 
\ (see \cite[Lem. 1.2(1.d)]{KrPaSo01})  
applied to $f:= \Res_{k \cA} $ and 
$g:= \init_\prec (\Res_{k \, \cB} ) $. 

\smallskip 

Finally we set $\cB:= b + d \, [0,1]^n \subset \R^n$ 
where $[0,1]$ denotes the unit interval of $\R$, 
for some 
$b \in \Z^n$ and $d \in \N$ such that $\cA_0, \dots, \cA_n 
\subset b + d \, [0,1]$. 
Then 
\,$ n(k) 
=\log \Big( \# ( k \, b + k\, d \, [0,1]^n \, \cap \Z^n) \Big) 
= \log  (k \, d+1)^n 
= O_k(\log k )$  
\ (here the notation $O_k$ refers to the dependence on
$k$) and so 
$$
h (\Res_{k \cA} ) = O_k( k^n \log k).  
$$
Alternatively, we could have obtained this from 
the inequality (\ref{c-e}) in the 
appendix. 

\smallskip

Toghether with the inequality (\ref{prod}) this implies that 
$$
\log \Big|\Res_\cA(f_0 , \dots, f_n) \Big | 
\, \le \,  
\sum_{i=0}^n MV_i (\cA) \, 
 \, \log ||f_i||_1 + 
O_k \bigg(\frac{\log k }{k}  \bigg), 
$$
from where we conclude by letting \ $ k \to \infty$. 
\end{Proof}

\bigskip

Let us consider some examples. 
For short we set 
$H(\cA):= H(\Res_\cA) $ and  $E(d) := 
\prod_{i=0}^n 
(\# \cA_i)^{MV_i(\cA)}$, and  
we also set 
$$
q(\cA):= 
\frac{\log E(\cA)}{ \log H(\cA)} 
$$ 
for the quotient
between the height of the resultant
and the estimate from  Theorem~\ref{h}.

\begin{exmpl} 
{\it Sylvester resultants.} 
{\rm For $d \in \N$ we let 
$$
\cA_0(d) = \cA_1(d):= \{ 0, 1, 2, \dots, d\} \ \subset \Z. 
$$
The corresponding resultant coincides with the 
Sylvester resultant   of two univariate polynomials of the same  
degree $d$. 
In this case 
 $MV_0(d)=MV_1(d) = d$ and   
$ \# \cA_0(d) = \#\cA_1(d) = d+1$, and so $E (d):= 
E\Big(\cA_0(d), 
\cA_1(d)\Big) 
= (d+1)^{2 \, d} $. 

\smallskip 

We compute the height $H(d):= H\Big(\cA_0(d), \cA_1(d)\Big) $ 
for $2 \le d \le 7$
with
the aid of Maple, and  
we collect the results in the following comparative table: 

\bigskip

\begin{center} 
{\small

\begin{tabular}{c|ccccccc} 
 \hline
&&&&&&& \\[-3.2mm]
\hspace{3mm} $d$  
\hspace{3mm} & \ &  $2$ &
$3$  
& $4$ & $5$ & $6$ & $7$ \\[1mm] 
\hline 
&&&&&&& \\[-2.7mm]
$H(d)$ & \ & $2$ & $3$ & $10$ & $23$ & $78$ & $274$ \\[1.5mm] 
$E(d)$  & \ & $81$ & $4,\hspace{-0.5mm}096$ & 
$390,\hspace{-0.5mm}625$ & $60,\hspace{-0.5mm}466,\hspace{-0.5mm}176$ & 
$13,\hspace{-0.5mm}841,\hspace{-0.5mm}287,\hspace{-0.5mm}201$ & 
$4,\hspace{-0.5mm}398,\hspace{-0.5mm}046,\hspace{-0.5mm}511,\hspace{-0.7mm}104$ \\[1.5mm] 
$q(d)$ &  \ &  $6.33$ & $7.57$ & $5.59$
& $5.71$ & $5.35$ & $5.18$ \\[1mm]  
\hline 

\end{tabular} 

}
\end{center}

} 
\end{exmpl} 

\smallskip

\begin{exmpl} \label{emiris} 
{\rm
We take this example from \cite[Exmpl. 3.5]{EmMo99}. 
Let 
\begin{eqnarray*} 
\cA_0 & := & \Big\{ (0,0), (1,1), (2,1), (1,0)\Big\},  \\[1mm] 
\cA_1 & := & \Big\{ (0,1), (2,2), (2,1), (1,0)\Big\},  \\[1mm] 
\cA_2 & := & \Big\{ (0,0), (0,1), (1,1), (1,0)\Big\}.    
\end{eqnarray*} 
Then $MV_0= 4$, $MV_1= 3$ and  $MV_2= 4$, so that 
$ E(\cA) = 4^4\, 4^3 \, 4^4$. 
On the other hand, we can compute 
the resultant using its expression in \cite[Exmpl. 3.19]{EmMo99}
as a quotient of determinants, 
and we obtain that 
$H(\cA) =8 $. 
Hence 
$$
H(\cA) = 8 \quad \quad , \quad \quad 
 E(\cA)= 4,\hspace{-0.8mm}194,\hspace{-0.6mm}304 
\quad \quad , \quad \quad 
q(\cA) = 7.33 \ . 
$$
For reference, the straightforward estimation 
of the Canny-Emiris formula gives (see the appendix below):  
$$
H(\cA) \le 4^{41} =
4,\hspace{-0.6mm}835,\hspace{-0.6mm}703,\hspace{-0.6mm}278,\hspace{-0.6mm}458,\hspace{-0.6mm}516,\hspace{-0.6mm}698,\hspace{-0.6mm}824,\hspace{-0.6mm}704
\ .  
$$

} 
\end{exmpl}

\begin{exmpl} 
{\rm
We take 
this example  from \cite[Exmpl. 2.1]{Sturmfels94}. 
Let 
\begin{eqnarray*} 
\cA_0 & := & \Big\{ (0,0), (2,2), (1,3)\Big\},  \\[1mm] 
\cA_1 & := &  \Big\{ (0,1), (2,0), (1,2) \Big\},  \\[1mm] 
\cA_2 & := & \Big\{ (3,0), (1,1)\Big\}.    
\end{eqnarray*} 
Then 
$MV_0= 5$, $MV_1= 7$ and  $MV_2= 7$, so that 
$ E(\cA) = 3^5\, 3^7 \, 2^7$. 
{}From the explicit monomial expansion of the resultant 
(see \cite[Exmpl. 2.1]{Sturmfels94}) we find that 
$H(\cA) = 14$ and so 
$$
H(\cA) = 14 \quad \quad , \quad \quad 
 E(\cA)= 68,\hspace{-0.6mm}024,\hspace{-0.6mm}448 \quad \quad , \quad \quad 
q(\cA) = 6.83 \ .
$$

} 
\end{exmpl}

\smallskip

These examples show that there is still some room for improvement 
over 
Theorem~\ref{h}. 
It is however possible that our estimate is quite sharp anyway: 
in spite of the large difference between $H(\cA)$ and 
$E(\cA)$ in the computed examples, 
the quotient $q(\cA)$ is quite small, and moreover it does not seem to
grow when $E(\cA) \to \infty$. 

In any case, it 
would be very interesting to have an {\em exact} expression for
$h(\Res_\cA)$ 
--- as remarked to me by B. Sturmfels --- or at least a non trivial
lower bound. 
Note that the only information that we dispose about the exact
value of the coefficients of $\Res_\cA$ is for the extremal ones,
which are equal to $\pm 1$ \cite[Cor. 3.1]{Sturmfels94}.

\begin{rem} 
While  a first version of this paper was circulating, C. D'Andrea 
(personal communication) 
obtained a non trivial lower bound for the height of the Sylvester 
resultant, and an improvement of the upper bound to 
$H\Big(\Res(f,g)\Big)  \le \max\{ \deg(f)!, \deg(g)!\}$. 

\end{rem} 

\medskip 

A final remark with respect to the Mahler measure.   
Let $X_{\cA} \subset \P^{m_0-1} \times \cdots \times \P^{m_n-1} $ 
be the projective toric variety associated to $\cA$ and
let
$\cL_i $ denote the equivariant line bundle which corresponds to
$\cA_i$
for $i=0, \dots, n$
as defined  in \cite[Section 8.1]{GeKaZe94}. 
Then 
the $\cA$-resultant coincides with 
the $(\cL_0, \dots, \cL_n)$-resultant of $X_\cA$, see 
\cite[Ch. 8, Prop. 1.5]{GeKaZe94}. 

In the context of Arakelov geometry, it is possible to equip these lines bundles  
with a canonical Hermitian  
metric following  \cite[Section 3.3]{Maillot00}. 
Then 
it is then natural to conjecture that 
$$
m\Big(\Res_{\cA_0, \dots, \cA_n} \Big) = h_{\overline{\cL}_0, \dots, 
\overline{\cL}_n} (X_\cA), 
$$
where $h_{\overline{\cL}_0, \dots, 
\overline{\cL}_n} (X_\cA)$ 
denotes the 
multiheight of  $X_\cA$ with respect to  the Hermitian 
line
bundles 
$ \overline{\cL}_0, \dots, 
\overline{\cL}_n $. 
This is suggested by an analogous result for the height of a
projective variety with respect to the Fubini-Study metric, 
see e.g. \cite[Thm. 3]{Soule92} or \cite[Cor. 2.4]{Philippon91}. 

If this is the case, together with  \cite[Prop. 7.11]{Maillot00} 
this would improve
Theorem~\ref{m} above
to  $m(\Res_{\cA_0, \dots, \cA_n})=0$.

\smallskip 

\section*{Appendix: Estimation of the height via the Canny-Emiris formula}

For purpose of easy reference, we 
establish herein 
the estimate for $h(\Res_{\cA})$ which follows from the
Canny-Emiris formula and the standard estimates for the behavior of 
the height of polynomials under addition, multiplication and
division.

\medskip

Assume that $L_\cA = \Z^n$ and set $Q:= \sum_{i=0}^n Q_i \subset
\R^n$. 
Given 
a coherent mixed subdivision of $Q$ and 
 a set
$$
\cE:= ( Q  + \delta ) \cap \Z^n,  
$$
where 
$\delta \in \Q^n$ is
sufficiently small and generic, there is a family of 
Canny-Emiris (square non singular) matrices $\cM_0, \dots, \cM_n$. 

For $j=0, \dots, n$, 
the given subdivision of $Q$ 
splits the set $\cE$ into a disjoint union 
$\cE= \cE_0(j) \cup \cdots \cup \cE_n(j)$.  
The elements in $\cE$ are in bijection with the rows of $\cM_j$, 
and to each 
 $ p\in \cE_i(j)$  corresponds a row of $\cM_j$ 
with 
exactly $m_i$ non zero entries, which consist
of the variables in  $U_i := \{ U_{i \, a} \, ; \, a \in \cA_i\}$. 
We refer to \cite[Sec. 7.6]{CoLiOS98} 
for the precise construction. 

Set $D_j:= \det(\cM_j)\in \Z[U_0, \dots, U_n] \setminus \{ 0\}$. 
The {\em Canny-Emiris formula} \cite[Ch. 7, Thm. 6.12]{CoLiOS98} states
that 
$ \Res_\cA = \gcd(D_0,\dots, D_n)$.

\smallskip

Then    
 $D_0$ 
is a multihomogeneous
polynomial of degree $N_i:= \# \cE_i (0) $ in each group of 
variables $U_i$ and 
of height bounded by $\sum_{i=0}^n N_i \, \log (m_i)$.
We have that 
\,$\Res_\cA | D_0$\, and so 
\,$m(\Res_\cA) \le m(D_0)$, which combined with 
\cite[Lem. 1.1]{KrPaSo01} gives 
\begin{equation} \label{c-e} 
h(\Res_\cA)  \, \le \,      
h(D_0)  +  \sum_{i=0}^n \Big( N_i + MV_i(\cA) \Big) \, \log (m_i) 
\, \le \, 
 \sum_{i=0}^n \Big( 2\, N_i + MV_i(\cA) \Big) \, \log (m_i).  
\end{equation}

Applied to Example~\ref{emiris},  this gives the stated estimate: 
then  $N_0 = N_1 = 4$ and  $N_2= 7$ (see
\cite[Exmpl. 3.5]{EmMo99}) and so the previous estimate gives 
 $H(\cA) \le 4^{2\cdot 15 + 11}= 4^{41} $. 

\smallskip 

In general, 
the estimate so obtained  is {\em much} worse than that of Theorem~\ref{h}, 
especially for $n\gg0$.
Consider e.g. 
$
\cA_i := \{0,\dots, d\}^n \subset \Z^n $ for 
$ i=0, \dots, n$.  
Then it is easy to show that 
Inequality (\ref{c-e}) gives  
$$
h(\Res_\cA) \le 
\Big(2\, ((n+1)\, d)^n + (n+1)\, d^n \Big)\, \log(d+1)
$$  
while Theorem~\ref{h} gives $h(\Res_\cA) \le (n+1) \, d^n \, \log(d+1)$.


\medskip 

\subsection*{Acknowledgements} I thank  Carlos D'Andrea 
for several helpful discu\-ssions, and for
having  communicated me his work on 
the height of the Sylvester resultant.


\end{document}